\documentclass{amsart}

\usepackage{times}
\usepackage[latin1]{inputenc}
\usepackage{graphicx}
\usepackage{amssymb}
\def\C{\mathbb{C}}
\def\R{\mathbb{R}}

\def\ds{\displaystyle}

\def\O{\mathcal{O}}

\newcommand{\Img}[2]{\includegraphics[width=#1truecm]{#2}}

\newtheorem{theorem}{Theorem}
\newtheorem{lemma}{Lemma}
\newtheorem{corollary}{Corollary}

\begin{document}

\date{November 11, 2006. Revised: July 11, 2007}
\title{On Differences of Zeta Values}
\author{Philippe Flajolet and Linas Vepstas}
\thanks{\emph{Address}: Algorithms Project, INRIA-Rocquencourt, F-78153 Le Chesnay, France}
\thanks{\emph{Corresponding author}: {\tt Philippe.Flajolet@inria.fr}}
\thanks{\emph{2000 Mathematics Subject Classification}:
11M06, 30B50, 39A05, 41A60}

\begin{abstract}
Finite differences of values of the Riemann zeta function at the integers
are explored. Such quantities, which occur as coefficients in Newton series representations, 
have surfaced in works of Bombieri--Lagarias, Ma{\'s}lanka, Coffey, B{\'a}ez-Duarte, Voros and others.
We apply the theory of N\"orlund-Rice integrals in conjunction with
the saddle-point method  and derive precise 
asymptotic estimates. The method extends to Dirichlet $L$-functions
and our estimates appear to be partly related to
earlier investigations surrounding Li's criterion for the Riemann hypothesis.
\end{abstract}

\keywords{Riemann zeta function, finite differences, 
asymptotic analysis, saddle point method, 
Li's criterion}

\maketitle

\section*{Introduction}

In recent  times, a variety of  authors have, for a variety  of reasons,
been led  to considering   properties  of  representations of   the
Riemann zeta function $\zeta(s)=\sum 1/n^s$  as a \emph{Newton interpolation
series}.  Amongst the many  possible  forms, we single out
the one relative to a regularized version of Riemann zeta, namely,
\begin{equation}\label{eq1}
\zeta(s)-\frac{1}{s-1}=\sum_{n=0}^\infty
(-1)^n b_n \binom{s}{n},
\end{equation}
where $\binom{s}{n}$ is a binomial coefficient:
\[
\binom{s}{n}:=\frac{s(s-1)\cdots(s-n+1)}{n!}.
\]
Corollary~\ref{newton-cor} in Section~\ref{conv-sec}
establishes that the representation~\eqref{eq1} is
valid throughout the complex plane,
its coefficients being determined by a general formula in the calculus of finite
differences (see~\cite{Jordan65,Milne81,Norlund54} and Section~\ref{newtser-sec} below):
\begin{equation}\label{eq2}
b_n=
n(1-\gamma-H_{n-1})-\frac12+\sum_{k=2}^n \binom{n}{k}
(-1)^k\zeta(k),
\end{equation}
Here,
$H_n=1+\frac12+\cdots\frac1n$ is a harmonic number.
Although the terms in the sum defining~$b_n$ become
exponentially large (of order close to~$2^n$), the values of the $b_n$ 
turn out to be exponentially small,
while exhibiting a curious oscillatory behavior. 
We shall indeed prove the estimate (Theorem~\ref{zetacoeff-thm} in Section~\ref{sadzeta-sec}) 
\begin{equation}\label{asympform}
b_n =\left(\frac{2n}{\pi}\right)^{1/4}e^{-2\sqrt{\pi n}}\cos\left(2\sqrt{\pi n}-\frac{5\pi}{8}\right)
+\O\left(n^{-1/4} e^{-2\sqrt{\pi n}}\right).
\end{equation}


Our \emph{first} motivation for investigating~\eqref{eq1}
and~\eqref{eq2} was an attempt by one of us, Linas
[2003, unpublished; available at
\verb|http://linas.org/math/poch-zeta.pdf|], to obtain
alternative and tractable expressions for the Gauss-Kuzmin-Wirsing operator of 
continued fraction theory 
In particular, Linas' computations at that time revealed that 
the $b_n$ tend rather fast to~0 and exhibit a surprising oscillatory pattern,
both facts crying for explanation. The present paper, essentially elaborated in early 2006,
represents the account of our joint attempts at understanding what goes on.

\emph{Second}, the zeta function has received attention in physics,
for its role in regularization and renormalization in 
quantum field theory. 
Motivated by such connections, 
Ma\'slanka introduced in~\cite{Maslanka01} 
what amounts to the Newton series representation of 
a regularized verson of $\zeta(s)$, namely, $(1-2s)\zeta(2s)$.
(Further numerical observations relative to the corresponding coefficients are
presented by this author in~\cite{Maslanka04}.)
The growth of coefficients
in Ma\'slanka's expansion has been subsequently investigated  by B\'aez-Duarte~\cite{Baez03}.
In particular, B\'aez-Duarte's analysis implies that the coefficients
decrease to~0
faster than any power of $1/n$. The methods we develop to derive Equation~\eqref{asympform}
can be easily adapted to yield refinements of the estimates of Baez-Duarte and Ma\'slanka.

A \emph{third} reason for  interest in the representation~\eqref{eq1} and the
companion           coefficients~\eqref{eq2}         is     \emph{Li's
criterion}~\cite{Li97} for the   Riemann Hypothesis (RH).   Let $\rho$
range over the  nontrivial zeros of~$\zeta(s)$.   Li's theorem asserts
that RH is true if and only if all members of the sequence
\[
\lambda_n=\sum_\rho \left[1-\left(1-\frac{1}{\rho}\right)^n\right], \qquad n\ge1,
\]
are nonnegative. Bombieri and Lagarias~\cite{BoLa99}
 offer an insightful discussion of Li's criterion. 
Coffey~\cite{Coffey05}, following Bombieri and Lagarias~\cite[Th.~2]{BoLa99}, has expressed the $\lambda_n$
as a  sum of  two terms, one  of which
 is an elementary variant $b_n$ defined by
\[
a_n(1,2):= \sum_{k=2}^n \binom{n}{k}(-1)^k (1-2^{-k})\zeta(k)-nH_{n-1}-\frac{n}{2}(\gamma-1+\log2).
\]
Theorem~2 of~\cite{Coffey05}  amounts to the property that the coefficients $a_n(1,2)$ 
decrease to~0. As we shall see in Section~\ref{dir-sec} and right after Equation~\eqref{coffeyS1}
in Section~\ref{concl-sec},
the methods originally developed for estimating~$b_n$
yield precise asymptotic information on~$a_n(1,2)$ as well.
Though the sums  we deal with 
count amongst the  far easier ones,
our precise asymptotic estimates may contribute to bring some clarity in this range of 
problems.

\smallskip

In this essay,  we  approach the problem of  asymptotically estimating
differences of  zeta  values by means of   a combination  of  two well
established   techniques.   We   start   from   a   contour   integral
representation of these  differences as defined by~\eqref{eq2}  (for this
technique, see   especially N\"orlund's treatise~\cite{Norlund54}  and
the study~\cite{FlSe95}),  then proceed to estimate the corresponding
complex  integral by  means of the classical  saddle-point method of asymptotic
analysis~\cite{deBruijn81,Olver74}. Our approach parallels a recent  paper
of Voros~\cite{Voros06} (motivated by Li's criterion), 
which our results supplement by providing a fairly 
detailed asymptotic analysis of differences of zeta values.

Section~\ref{newtser-sec} presents the construction of a Newton series for the zeta function,
and presents generating functions for its coefficients. This is followed,  in Section~\ref{exper-sec},
by a brief examination of numerical results. Section~\ref{norl-sec}
gives the N\"orlund integral representation for the coefficients,
of which  Section~\ref{sadzeta-sec} provides a careful saddle-point analysis.
 The convergence of the Newton series
representation~\eqref{eq1} is then discussed in Section~\ref{conv-sec}, and
Section~\ref{dir-sec} develops the corresponding analysis
for Dirichlet $L$-functions.
 We end with a conclusion, Section~\ref{concl-sec}, outlining other applications of N\"orlund integrals
in the realm of finite differences and zeta functions.

\section{Newton series and zeta values}\label{newtser-sec}

This section defines the Newton series for the Riemann zeta
that is to be studied, demonstrates some of its basic properties,
and gives some generating functions for its coefficients.
In this paper, a Newton series is defined as

\begin{equation}\label{gennew}
\Phi(s)=\sum_{n=0}^\infty (-1)^n c_n \binom{s}{n}.
\end{equation}
Given a function $\phi(s)$, one may attempt to represent it 
in some region of the complex plane by means of 
such a series. Since the series $\Phi(s)$ terminates at $s=0,1,2,\ldots$,
the conditions $\phi(m)=\Phi(m)$ at the nonnegative integers
imply that the candidate sequence $\{c_n\}$ is linearly related to the
sequence of values $\{\phi(m)\}$ by
\[
\phi(m)=\sum_{n=0}^m (-1)^n c_n \binom{m}{n}.
\]
The triangular system can then be inverted to 
give (by the binomial transform~\cite[p.~192]{GrKnPa89},
or its Euler transform version~\cite[p.~43]{Riordan68},
or by direct elimination)
\begin{equation}\label{gencoeff}
c_n = \sum_{k=0}^{n} \binom{n}{k} (-1)^k \phi(k), \qquad n=0,1,2,\ldots\,.
\end{equation}
This choice of coefficients for~\eqref{gennew}
determines the Newton series \emph{associated} to $\phi$.
The coincidence of 
the function $\phi$ and its associated series
$\Phi$ is, by construction, granted at least
at all the nonnegative integers. 
The validity
of $\Phi(s)=\phi(s)$  is often found to extend to large
parts of the complex plane, but this fact requires specific
properties  beyond the mere convergence of the series in~\eqref{gennew}.


In the case of the Newton series for $\zeta(s)-1/(s-1)$, the general
relation~\eqref{gencoeff} provides the coefficients in the form
\begin{equation}\label{bn0}
b_n = s_0-ns_1+\sum_{k=2}^n \binom{n}{k} (-1)^k \left[\zeta(k)-\frac{1}{k-1}
\right],
\end{equation}
where
\begin{equation}\label{bn1}
s_0=\left[\zeta(s)-\frac{1}{s-1}\right]_{s=0}=\frac12,
\qquad
s_1=\lim_{s\to 1} \left[\zeta(s)-\frac{1}{s-1}\right] = \gamma.
\end{equation}
The harmonic numbers appear as\footnote{%
	Identity~\eqref{bn2} is easily deduced from the classical partial fraction
	decomposition (see, e.g., \cite[p.~188]{GrKnPa89}),
\[
\frac{n!}{x(x+1)\cdots(x+n)}=\sum_k \binom{n}{k} \frac{(-1)^k}{x+k},
\]
upon letting $x\to-1$.
}
\begin{equation}\label{bn2}
\sum_{k=2}^n \binom{n}{k}\frac{(-1)^k}{k-1}=1-n+nH_{n-1}.
\end{equation}
Equations \eqref{bn0}, \eqref{bn1},   \eqref{bn2} 
then entail that   the
$b_n$, as  defined by~\eqref{eq2}, are indeed the  coefficients  of the Newton series
associated to~$\zeta(s)-1/(s-1)$. A proof that the equality
$\Phi(s)=\zeta(s)-1/(s-1)$ holds for all complex~$s$ is given in Section~\ref{conv-sec},
following the asymptotic analysis of the coefficients~$b_n$ and based on
a theorem of Carlson.

\smallskip

Before engaging in a detailed study of the~$b_n$, we note a few 
simple facts regarding their elementary properties. Consider the quantities
\begin{equation}\label{deldef}
\delta_n:=\sum_{k=2}^n \binom{n}{k} (-1)^k \zeta(k),
\end{equation}
which represent the nontrivial sums in the definition~\eqref{eq2} of~$b_n$
and are, up to minor adjustments, differences of zeta values at the
integers. 
Expanding the zeta function according to its definition and 
exchanging the order of summations in the resulting double sum yields
\begin{equation}\label{simpdel}
\delta_n = \sum_{\ell\ge1} \left[\left(1-\frac{1}{\ell}\right)^n-1+\frac{n}{\ell}\right].
\end{equation}
This rather simple sum shows a remarkably complex behavior;
elucidating its behavior is one of the principal topics of this paper.

The ordinary generating function for the sequence $\{\delta_n\}$ 
is also of interest.
Given the classical expansion~\cite[{\S6.3}]{AbSt73} 
of the logarithmic derivative $\psi(z)=
\Gamma'(z)/\Gamma(z)$ of the Gamma function,
\[
\psi(1+z)+\gamma=\zeta(2)z-\zeta(3)z^2+\zeta(4)z^3-\cdots,
\]
one finds, by the usual generating function translation of the 
Euler transform~\cite[p.~311]{Milne81} or by an immediate verification
based on the binomial theorem:
\begin{equation}\label{delogf}
\sum_{n\ge2}\delta_n z^n =\frac{z}{(1-z)^2}\left[\psi\left(\frac{1}{1-z}\right)
+\gamma\right].
\end{equation}
The exponential generating function for the sequence  $\{\delta_n\}$ reflects~\eqref{simpdel} and is even simpler:
\begin{equation}\label{delegf}
\sum_{n\ge2}\delta_n \frac{z^n}{n!}
=e^z\sum_{n\ge2} \zeta(n) \frac{(-z)^n}{n!}=e^z \sum_{\ell\ge1}
\left[e^{-z/\ell}-1+\frac{z}{\ell}\right] .
\end{equation}
(The first equality is easily verified by expanding $e^{z}$ and expressing the coefficient
of~$z^n/n!$ in the product 
as a convolution, itself seen to coincide with~\eqref{deldef}.)

\section{Experimental analysis}\label{exper-sec}

Detailed experiments on the  coefficients $b_n$ conducted by one of us  are
at the origin of the present paper and we briefly discuss these 
since they illustrate some concrete numerical aspects of the sequence $(b_n)$
while being potentially useful for similar problems.
As it is usual when dealing with finite differences,
the alternating binomial sums giving the $b_n$ involve exponential
cancellation since the binomial coefficients get almost as large as $2^n$. 
We started by conducting evaluations of the $b_n$ up to $n\approx 5000$,
which requires determining zeta values up to several thousand digits of precision.
Note that the zeta values can be computed rapidly to extremely high
precision using several efficient algorithms (e.g.,~\cite{CoViZa00}), which are available in 
symbolic computation packages ({\sc Maple}) and numerical libraries
({\sc Pari/Gp}).

A quick inspection of numerical data immediately reveals two features
of the constants~$b_n$: they are oscillatory with a slowly increasing
(pseudo)period and their
absolute values are rapidly decreasing. For instance\footnote{%
	The notation $x\doteq y$ designates a numerical approximation of
	$x$ by $y$ to the last decimal digit stated. 
}:
\[\renewcommand{\arraycolsep}{2truept}
\begin{array}{lllllllll}
b_1&\doteq&     -7.72156\cdot 10^{-2},\quad &
b_2&\doteq&	-9.49726\cdot 10^{-3},\quad &
b_5&\doteq& 	 +7.15059\cdot 10^{-4}, \\
b_{10}&\doteq&  -2.83697\cdot 10^{-5}, \quad &
b_{20}&\doteq& 	 +2.15965\cdot 10^{-9},\quad &
b_{50}&\doteq&   -1.08802\cdot 10^{-11}.
\end{array}\]

A numeric fit of the oscillatory behavior of the function was first made.
There are sign changes in the sequence $\{b_n\}$ at
\[
n=3,7,13,21,29,40,52,65,80,97,115,135,157,180,\ldots\,,
\]
the values growing roughly quadratically. A good fit for the $k$th sign-change 
was found experimentally to be of the form $q(k)= \frac{\pi}{4}k^2+\O(k)$.
The quadratic polynomial is then easily 
inverted to give an approximate oscillatory behavior of the $b_n$.
%
%
Once the oscillatory behavior had been disposed of, the task of quantifying
the general trend in the overall decrease of the sequence became easier. 
These observations then led us to conjecture
\begin{equation}\label{conj}
\beta(n):=  \cos\pi\left(2\sqrt{\frac{n}{\pi}}+L\right)  e^{-K\sqrt{n}},
\qquad K=3.6\pm 0.1,
\end{equation}
as  a   rough approximation  to~$b_n$, where $L$ is some real constant.
This numerical fit then greatly helped us find the
main estimate~\eqref{asympform} above, which has
the exact value $K=\sqrt{2\pi}\doteq 3.54490$ and
an additional $n^{1/4}$ factor modulating the exponential.


\section{The N\"orlund integral representation}\label{norl-sec}

Our approach to the asymptotic estimation of the $b_n$ relies on
a complex integral representation of finite differences of an analytic function,
to be found in N\"orlund's classic treatise~\cite[\S VIII.5]{Norlund54} first published in 1924.
In computer science, this representation was popularized by Knuth~\cite[p.~138]{Knuth98a},
who attributed it to S.O.~Rice, so that
it also came to be known as ``Rice's method''; see~\cite{FlSe95}
for  a review. 

\begin{lemma}\label{norlem}
 Let $\phi(s)$ be holomorphic in the half-plane 
$\Re(s)\ge n_0-\frac12$. Then the finite differences of the sequence
$(\phi(k))$ admit the integral representation
\begin{equation}\label{norint}
\sum_{k=n_0}^n \binom{n}{k}(-1)^k\phi(k)=\frac{(-1)^{n}}{2\pi i}\int_C 
\phi(s) \frac{n!}{s(s-1)\cdots(s-n)}\, ds,
\end{equation}
where the contour of integration~$C$ encircles the integers $\{n_0,\ldots,n\}$
in a positive direction and is contained in~$\Re(s)\ge n_0-\frac12$.
\end{lemma}
\begin{proof}
The integral on the right of~\eqref{norint} is the sum of its residues at~$s=n_0,\ldots,n$,
which precisely equals the sum on the left.
\end{proof}

An immediate consequence is the following representation for the 
differences of zeta values ($\delta_n$ is as in~\eqref{deldef}):
\begin{equation}\label{norbn}
\delta_n\equiv \sum_{k=2}^n \binom{n}{k}(-1)^k\zeta(k)
=\frac{(-1)^{n-1}}{2\pi i}\int_{3/2-i\infty}^{3/2+i\infty} 
\zeta(s) \frac{n!}{s(s-1)\cdots(s-n)}\, ds. 
\end{equation}
(Choose $T>n$ and consider the finite contour (negatively oriented) consisting of the line from
$\frac32-iT$ to $\frac32+iT$, followed by the clockwise arc of $|s|=\sqrt{T^2+9/4}$ that
lies to the right of the given line. The contribution of the circular arc is $\O(T^{-n})$ 
and thus vanishes as $T\to+\infty$.)

Since~$b_n$ is $\delta_n$ plus a correction term (see Equation~\eqref{eq2})
and $\delta_n$ admits the integral representation~\eqref{deldef},
the first step of our analysis is to move the line of integration further to the left.
It is well known that the Riemann zeta function is of finite order 
in any right half-plane~\cite[\S5.1]{Titchmarsh86},
that is, $|\zeta(s)|=\O(|s|^A)$ uniformly as $|s|\to\infty$, for some~$A$ depending on the
half-plane under consideration.
As a consequence, the integral of~\eqref{norbn} remains 
convergent, when taken along any vertical line left of~0,
as soon as~$n$ is large enough. Under these conditions, it is possible
to replace the line of integration $\Re(s)=\frac32$ by the line~$\Re(s)=-\frac12$, upon
taking into account the residues of a double pole at $s=1$ and a simple
pole at $s=0$.
We find in this way
\[
\delta_n=(-1)^{n-1} (R_1+R_0)+\frac{(-1)^{n-1}}{2\pi i}\int_{-1/2-i\infty}^{-1/2+i\infty} 
\zeta(s) \frac{n!}{s(s-1)\cdots(s-n)}\, ds, \]
where, as shown by a routine calculation:
\[
(-1)^n R_0=-\frac12,
\qquad
(-1)^nR_1=-n(1-\gamma-H_{n-1}).
\]
The residues thus compensate exactly 
for the difference between~$\delta_n$ and~$b_n$, so that
\begin{equation}\label{norbn2}
b_n=
\frac{(-1)^{n-1}}{2\pi i}\int_{-1/2-i\infty}^{-1/2+i\infty} 
\zeta(s) \frac{n!}{s(s-1)\cdots(s-n)}\, ds.
\end{equation}

\section{Saddle-point analysis of zeta differences}\label{sadzeta-sec}

The integrand in~\eqref{norbn2} is free of 
singularities on the left-hand side, and is thus a good candidate for estimation
by the saddle-point (or steepest-descent) method.  This evaluation 
is the main topic of this section, and culminates with the derivation of
one of the principal results of this paper, 
Theorem~\ref{zetacoeff-thm}.

We make use of the functional equation of
the Riemann zeta function under the form
\begin{equation}\label{funzeta}
\zeta(s) =2\Gamma(1-s)(2\pi)^{s-1}\sin \frac{\pi s}{2}\zeta(1-s).
\end{equation}
From this relation, upon performing the change of variables $s\mapsto -s$,
we obtain
\begin{equation}\label{mainint}
b_n = -\frac{1}{\pi i} 
\int_{1/2-i\infty}^{1/2+i\infty} 
(2\pi)^{-s-1} \sin\left(\frac{\pi s}{2}\right)\zeta(1+s)
\frac{n!\,\Gamma(1+s)}{s(s+1)\cdots(s+n)}\, ds,
\end{equation}
which is the starting point of our asymptotic analysis.

The integral representation~\eqref{mainint} has several noticeable features. 
First, the integrand has no singularity at all in $\Re(s)\ge\frac12$
and it appears to decay exponentially fast towards $\pm i \infty$
(see Equations~\eqref{zetapp}, \eqref{stirapp}, and~\eqref{stirapp} below).
This means that one can freely choose the abscissa $c$ (with $c\geq\frac12$)
in the representation 
\begin{equation}\label{mainint2}
b_n = -\frac{1}{\pi i} 
\int_{c-i\infty}^{c+i\infty} 
(2\pi)^{-s-1} \sin\left(\frac{\pi s}{2}\right)\zeta(1+s)
\frac{n!\, \Gamma(s+1)}{s(s+1)\cdots(s+n)}\, ds.
\end{equation}
The very absence of singularities calls for an application 
of the saddle- point method.

The factor $\zeta(1+s)$ remains bounded in modulus by a constant,
and is in fact barely distinguishable from~1, as $\Re(s)$ increases,
since
\begin{equation}\label{zetapp}
\zeta(s)=1+\O\left(2^{-\Re(s)}\right), \qquad \Re(s)\ge\frac32.
\end{equation}
Also, for large $|s|$, the complex version of Stirling's formula 
applies:
\begin{equation}\label{stirapp}
\Gamma(1+s) = s^s e^{-s}\sqrt{2\pi s}\left(1+\O(|s|^{-1})\right),\qquad \Re(s)\ge0.
\end{equation}
Finally, the sine factor increases exponentially along vertical lines:
one has
\begin{equation}\label{sinapp}
2i\sin\frac{\pi s}{2} =-\exp\left(-i\frac{\pi s}{2}\right)
+\O\left(e^{-\pi \Im(s)/2}\right), \qquad \Im(s)\ge0,
\end{equation}
with a conjugate approximation holding for $\Im(s)<0$.

In anticipation of applying saddle-point methods, 
the approximations~\eqref{zetapp}, \eqref{stirapp}, and~\eqref{sinapp}
then suggest the function $e^{\omega(s)}$
 as a simplified model of the integrand in the upper half-plane, where
\begin{equation}\label{omdef}
\omega(s)=-s\log(2\pi)-i\frac{\pi s}{2}+\log\frac{n! \Gamma(s)^2}{\Gamma(s+n)}.
\end{equation}
We shall demonstrate shortly that the location of the appropriate 
saddle-points in the complex plane scale as $\sqrt{n}$, which may be confirmed 
by numerical experiments. Therefore, in performing an asymptotic analysis,
it is appropriate to perform a change of variable $s=x\sqrt{n}$, 
and expand in descending powers of $n$, presuming $x$ to be approximately constant.
We find, uniformly for $x$ in any compact region of $\Re(x)>0$, $\Im(x)>0$:
\begin{equation}\label{om012}
\left\{\begin{array}{lll}
\omega(x\sqrt{n})&=&\ds x\sqrt{n}\left[2\log x-2-\log(2\pi)-\frac12i\pi \right]+\frac12\log n
\\
&&\ds \qquad -\log x +\log(2\pi)-\frac12x^2+\O(n^{-1/2})
\\
\omega'(x\sqrt{n})&=& \ds \left[-\log(2\pi) -i\frac{\pi}{2}+2\log x\right]
-\left(x+\frac1x\right)\frac{1}{\sqrt{n}}+\O(n^{-1}).
\\
\omega''(x\sqrt{n})&=&\ds \frac{2}{x\sqrt{n}}+\O(n^{-1}).
\end{array}\right. 
\end{equation}
(The symbolic manipulation system {\sc Maple} is a great help in such computations.)

From the second line of~\eqref{om012}, an approximate root of $\omega'(s)$
is obtained by choosing the particular value~$x_0$ of $x$ that cancels 
$\omega'$ to main asymptotic order:
\begin{equation}\label{x0}
x_0=e^{i\pi/4}\sqrt{2\pi}.
\end{equation}
This corresponds to the following value for $s$,
\begin{equation}\label{sad0}
\sigma\equiv \sigma(n)= x_0\sqrt{n} = 
(1+i)\sqrt{\pi n},
\end{equation}
which thus is also an approximate saddle-point for $e^{\omega(s)}$.
The substitution of this value given the first line of~\eqref{om012} then leads to
\begin{equation}\label{appsad}
\exp\left(\omega(\sigma(n))\right)
=\exp\left(2i\sqrt{\pi n}\right) \cdot \exp\left(-2\sqrt{\pi n}\right)\cdot \Pi(n),
\end{equation}
where $\Pi$ is an
unspecified factor of at most polynomial growth. 
By using a suitable contour that passes though~$\sigma(n)$,
we thus expect the quantity in~\eqref{appsad} to be an approximation
(up to polynomial factors again) of~$b_n$. This back-of-the-envelope calculation
does predict  the exponential decay of~$b_n$ as $\exp\left(-2\sqrt{\pi n}\right)$,
in a way consistent with numerical data, while 
the fluctuations,
$
\sin\left(2\sqrt{\pi n}+\O(1)\right),
$
are seen to be in stunning agreement with the empirically obtained formula~\eqref{conj}.

%

\begin{figure}

\begin{center}
\begin{footnotesize}
\setlength{\unitlength}{0.6truecm}
\begin{picture}(6,10)(-1,-5)
\thicklines
\put(-1,0){\line(1,0){6}}
\put(0,-5){\line(0,1){10}}
\put(0,0.2){~$0$}
\put(1,4.2){$\Re(s)=c_1\sqrt{n}$}
\put(3,0.2){$\Re(s)=c_2\sqrt{n}$}
\put(2,2){\circle*{0.2}}
\put(2,2){$\quad\ds \sigma=e^{i\pi/4}\sqrt{2\pi n}$}
\put(1,3){\line(1,-1){2}}
\put(1,3){\line(0,1){2}}
\put(3,1){\line(0,-1){1}}
\put(2,2){\vector(1,-1){0.5}}
\put(2,2){\vector(-1,1){0.5}}
\thinlines
\put(1,3){\line(0,-1){3}}
\thicklines
\put(2,-2){\circle*{0.2}}
\put(2,-2){$\quad\ds \overline\sigma=e^{-i\pi/4}\sqrt{2\pi n}$}
\put(1,-3){\line(1,1){2}}
\put(1,-3){\line(0,-1){2}}
\put(3,-1){\line(0,1){1}}
\put(2,-2){\vector(1,1){0.5}}
\put(2,-2){\vector(-1,-1){0.5}}
\thinlines
\put(1,-3){\line(0,1){3}}
\end{picture}
\end{footnotesize}
\qquad
\Img{7.5}{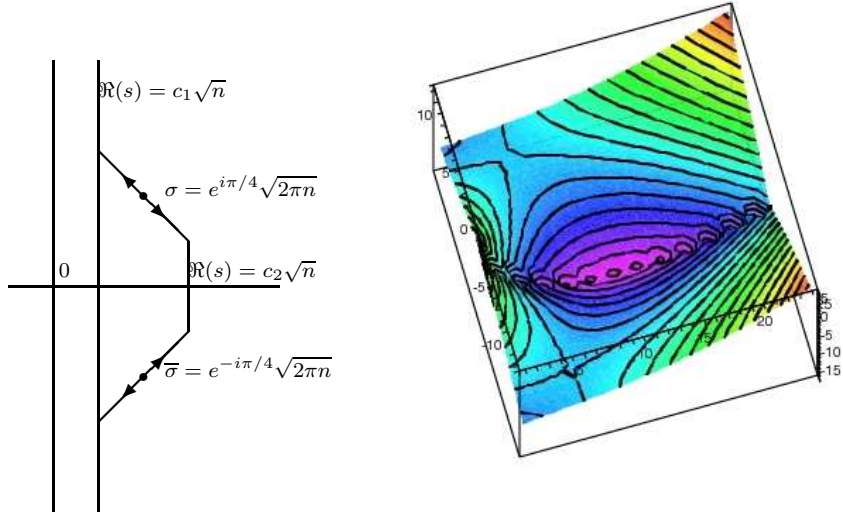}
\end{center}
\caption{\label{sad-fig}
Left: The saddle-point contour used for estimating $b_n$.
The arrows point at the directions of steepest descent
from the saddle-points.
Right: the landscape of the logarithm of the modulus of the integrand
in the representation of~$b_n$ for $n=10$.}

\end{figure}

\smallskip
We  must   now fix the   contour  of  integration  and  provide  final
approximations. The   contour  adopted  (Figure~\ref{sad-fig})    goes
through the saddle- point $\sigma=\sigma(n)$ and symmetrically  through
its  complex conjugate $\overline\sigma=\overline{\sigma(n)}$.  In the
upper half-plane, it traverses $\sigma(n)$ along  a line  of steepest
descent   whose   direction, as   determined    from the  argument  of
$\omega''(\sigma)$,  is at an    angle of $\frac{5\pi}{8}$    with the
horizontal  axis.  The contour   also  includes parts of two  vertical
lines of  respective abscissae $\Re(s)=c_1\sqrt{n}$ and~$c_2\sqrt{n}$,
where
\[
0<c_1<\sqrt{\pi}<c_2<2\sqrt{\pi}.
\] 

The  choice of the abscissae, $c_1$  and $c_2$, is
not critical (it is even possible to adapt the analysis to $c_1=c_2=\sqrt{\pi}$).
One verifies  easily, from crude approximations, that
the contributions arising from the vertical parts of the contour are
$\O(e^{-L_0\sqrt{n}})$,    for   some $L_0>2\sqrt{\pi}$,   i.e.,   they  are
exponentially small in the  scale of the problem:
\begin{equation}\label{vert}
\int_{\operatorname{vertical}} = \O\left(e^{-L_0\sqrt{n}}\right)\qquad L_0>2\sqrt{\pi}.
\end{equation}

The slanted part of
the contour is such that all the estimates of~\eqref{om012} apply. 
The scale of the problem is dictated by the value of $\omega''(\sigma)$,
 which is of order $\O(n^{-1/2})$. This indicates that the ``second order'' scaling 
to be adopted is~$n^{1/4}$. Accordingly, we set
\begin{equation}\label{cregion}
s=(1+i)\sqrt{\pi n}+e^{5i\pi/8}yn^{1/4}.
\end{equation}
Define the \emph{central region} of the slanted part of the contour by the condition
that $|y|\le \log^2 n$.
Upon slightly varying the value of~$x$ around~$x_0$, one verifies from~\eqref{om012} that,
for large~$n$, the quantity
\[
\Re\left(\frac{1}{\sqrt{n}}
\omega\left(x_0\sqrt{n}+e^{5i\pi/8}t\sqrt{n}\right)\right)
\]
is an upward concave function of~$t$ near~$t=0$. There results,
in the complement of the central part, $|y|\ge \log^2n$,
the approximation
\[
\left|\exp\left(\omega\left(x_0\sqrt{n}+e^{5i\pi/8}yn^{1/4}\right)\right)\right|
< e^{\omega(x_0\sqrt{n})}
\cdot \exp\left(-L_1 \log^2 n\right),\qquad L_1>0.
\]
Figuratively:
\begin{equation}\label{cent}
\int_{\operatorname{slanted}}=\int_{\operatorname{central}}+\O\left(\exp\left(-L_1 \log^2 n\right)\right).
\end{equation}
Thus, from~\eqref{vert}  and~\eqref{cent},  only the central   part of the
slanted region matters asymptotically. This applies to $e^{\omega(s)}$ but also
to the full integrand of the representation~\eqref{mainint2} of $b_n$,
given the approximations~\eqref{zetapp}--\eqref{om012}.

\smallskip

We are finally ready to reap the crop. Take the integral representation
of~\eqref{mainint2} with the contour deformed as indicated in Figure~\ref{sad-fig} and
let $b_n^{+}$ be the contribution arising from the upper half-plane, to
the effect that
\begin{equation}\label{conjug}
b_n=2\Re(b_n^+),
\end{equation}
by conjugacy. In the central region, 
\[
s=x_0\sqrt{n}+e^{5i\pi/8}yn^{1/4},
\]
the integrand of~\eqref{mainint2} becomes
\begin{equation}\label{phis}
\hbox{\small $\ds\left(-\frac{1}{\pi i}\right)\cdot 
(2\pi)^{-1}\cdot 
\left(-\frac{1}{2i}\right)\cdot
 \left(1+\O(2^{-\sqrt{\pi n}})\right)\cdot 
\frac{x_0}{\sqrt{n}} \cdot e^{\omega(x_0\sqrt{n})} \cdot
e^{-y^2/\sqrt{2\pi}}\left(1+\O\left(\frac{1}{\sqrt{n}}\right)\right)$}.
\end{equation}
The various factors found there (compare~\eqref{mainint2} to $e^{\omega(s)}$
with $\omega(s)$ defined in~\eqref{omdef})
are in sequence: the Cauchy integral prefactor; 
the correction $(2\pi)^{-1}$ to the functional equation of Riemann zeta;
the factor $-1/(2i)$ relating the sine to its
exponential approximation; the approximation of Riemann zeta;
the correction $s/(s+n)$ of the Gamma factors; 
the main term $e^{\omega(\sigma)}$;
the anticipated local Gaussian approximation;
the errors resulting from approximations~\eqref{zetapp}--\eqref{om012},
which are of relative order $\O(n^{-1/2})$. Upon completing the tails of the 
integral and neglecting exponentially small corrections,
we get
\begin{equation}\label{asy0}
 b_n^+ =K_0 e^{\omega(x_0\sqrt{n})}
\frac{x_0}{\sqrt{n}}
\int_{-\infty}^{+\infty} e^{-y^2/\sqrt{2\pi}}\,dy
\cdot\left(e^{5i\pi/8}n^{1/4}\right)\left(1+\mathcal{O}\left(\frac{1}{\sqrt{n}}\right)\right),
\end{equation}
where~$K_0=-1/(4\pi^2)$    is the constant  factor  of~\eqref{phis}, while the factor
following   the integral   translates the   change   of   variables:
$ds=e^{5i\pi/8}n^{1/4}dy$. 

The  asymptotic form of~$b_n$ is now completely
determined by~\eqref{conjug} and~\eqref{asy0}.
We have obtained:

\begin{figure}
\begin{center}
\Img{7}{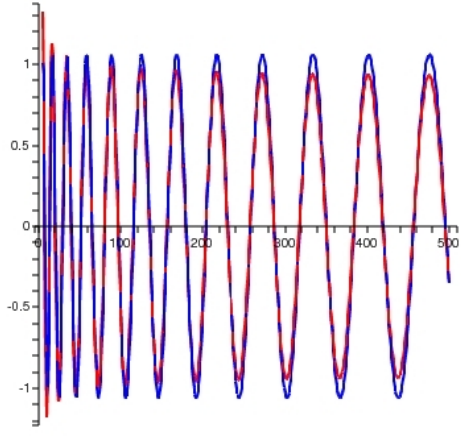}
\end{center}
\caption{\label{comparapp-fig} A comparative plot of $b_n$ and 
the main term of its approximation~\eqref{thmz},
both multiplied by $e^{2\sqrt{\pi n}}n^{-1/4}$, for $n=5\,.\,.\,500$.}
\end{figure}

\begin{theorem} \label{zetacoeff-thm}
The Newton coefficient $b_n$ of $\zeta(s)-1/(s-1)$ defined in~\eqref{eq2}
satisfies 
\begin{equation}\label{thmz}
b_n = \left( \frac{2n}{\pi}\right)^{1/4}
e^{-2\sqrt{\pi n}}\cos\left(2\sqrt{\pi n}-\frac{5\pi}{8}\right)
+\mathcal{O}\left(e^{-2\sqrt{\pi n}}n^{-1/4}\right).
\end{equation}
\end{theorem}

The agreement between asymptotic and exact values is quite good,
even for small values of~$n$ (Figure~\ref{comparapp-fig}).

\smallskip

In summary, the foregoing developments justify the validity of applying the saddle-point formula 
to the N\"orlund Rice
integral representation~\eqref{mainint2} of zeta value differences. Under
its general form, this formula (which, without further assumptions, remains a heuristic)
reads
\begin{equation}\label{sadfor}
\int e^{-Nf(x)} \, dx =\sqrt{\frac{2\pi}{N f''(x_0)}}
e^{-N f(x_0)}\left(1+\mathcal{O}\left(\frac{1}{N}\right)\right).
\end{equation}
Here the analytic function $f(x)$ must have a (simple) saddle-point at~$x_0$,
that is, $f'(x_0)=0$ and $f''(x_0)$ is the second derivative of~$f$ at the saddle-point.
In the case of differences of zeta values, 
the appropriate scaling parameter is $s=x\sqrt{n}$ corresponding to $N=\sqrt{n}$, and
the function~$f$ is
\[
f(x)=\lim_{n\to\infty} \frac{1}{\sqrt{n}}\omega\left(x\sqrt{n}\right),
\]
up to smaller order corrections that we could treat as constants in the range of the saddle-point.
As we shall see in Section~\ref{dir-sec}, this paradigm
adapts 
to sums involving Dirichlet $L$-functions.

\section{Convergence of the Newton series of zeta}\label{conv-sec}

The fact that the coefficients $b_n$ decay to zero faster than any polynomial in~$1/n$ 
implies that the Newton series
\begin{equation}\label{newtonz}
\Phi(s)=\sum_{n=0}^\infty (-1)^n b_n \binom{s}{n},
\end{equation}
with $b_n$ given by~\eqref{eq2},
converges throughout the complex plane, and consequently defines an entire function.
Set $Z(s):=\zeta(s)-1/(s-1)$ with~$Z(1)=\gamma$. We have, by construction $\Phi(s)=Z(s)$
at $s=0,1,2\ldots\,$, but the relation between $\Phi$ and $Z$ at other points is still unclear.

\begin{corollary}\label{newton-cor} The Newton series of~\eqref{newtonz} 
is a convergent representation of the function $\zeta(s)-1/(s-1)$
valid at all points $s\in\C$.
\end{corollary}
\begin{proof} 
Here  is our favorite  proof.  A   classic theorem  of Carlson  (for a
discussion        and      a     proof,     see,      e.g.,    Hardy's
Lectures~\cite[pp.~188-191]{Hardy78}          or          Titchmarsh's
treatise~\cite[\S5.81]{Titchmarsh39})  says the following:
\emph{Assume that $(i)$~$g(s)$ is
analytic and such that
\[
\left|g(s)\right|<C e^{A|s|},
\]
where $A<\pi$, in the right half-plane of complex values of~$s$,
and $(ii)$~$g(0)=g(1)=\ldots=0$. Then $g(s)$ vanishes identically.}

To complete the proof, it  suffices to apply  Carlson's theorem to the
difference    $g(s)=\Phi(s+2)-Z(s+2)$.
Condition~$(ii)$ is    satisfied  by   construction  of    the  Newton
series. Condition~$(i)$ results from  the fact that $Z(s+2)$ is $\mathcal{O}(1)$
while  a general bound  due to N\"orlund 
(Equation~(58) of~\cite[p.~228]{Norlund54})
and valid for all convergent Newton series
asserts that $|\Phi(s+2)|$  is of 
growth at most $e^{\frac{\pi}{2}|s|}$, throughout $\Re(s)>-\frac12$.
%
%
%
\end{proof}

An alternative proof can be given starting from a contour integral representation for the remainder of 
a general Newton series given~\cite[p.~223]{Norlund54}. 
Yet another proof derives from a turnkey theorem of N\"orlund, quoted
in~\cite[p.~311]{Milne81}:
\emph{In order that a function $F(x)$ should admit
a Newton series development, it is necessary and sufficient that $F(x)$
should be holomorphic in a certain half-plane $\Re(x)>\alpha$ and should there 
satisfy the inequality $|F(x)|<C 2^{|x|}$, where~$C$ is a fixed positive number.}
In a short note, B\'aez-Duarte~\cite{Baez03}
justified a similar looking Newton series representation of the zeta function due to Ma\'slanka---however
his bounds on the Newton coefficients are less precise than ours and his arguments (based on
a doubly indexed sequence of polynomials) seem to be somewhat problem-specific.

\section{Dirichlet $L$-functions}\label{dir-sec}

The methods employed to deal with differences of zeta values have a
more general scope, and we may reasonably expect them to be applicable to
other kinds of Dirichlet series.
Such is indeed the case
for any Dirichlet $L$-function,
\[
L(\chi,s)=\sum_{n=1}^\infty \frac{\chi(n)}{n^s},
\]
where $\chi$ is a multiplicative character 
of some period~$k$, that is,
for all integers~$m,n$, one has:
$\chi(n+k)=\chi(n)$,  $\chi(mn)=\chi(m)\chi(n)$, $\chi(1)=1$, and
$\chi(n)=0$ whenever~$\gcd(n,k)\not=1$.

%

Let $\zeta(s,q)$ be the Hurwitz
zeta function defined by
\begin{equation}\label{defhur}
\zeta(s,q)=\sum_{n=0}^{\infty}\frac{1}{(n+q)^{s}}\end{equation}
 Any Dirichlet $L$-function may be represented  as a combination of Hurwitz zeta functions,
\begin{equation}
L(\chi,s)=\frac{1}{k^{s}}\sum_{m=1}^{k}\chi(m)\zeta\left(s,\frac{m}{k}\right),
\label{eq:L-Hurwitz}\end{equation}
where $k$ is the period of $\chi$.
In particular, the coefficients of the Newton series  for $L(\chi,s)$ 
%
%
are simple linear combinations of the quantities
\begin{equation}\label{defA}
A_{n}(m,k)=\sum_{\ell=2}^{n}
\binom{n}{\ell}
(-1)^{\ell}\frac{\zeta\left(\ell,\frac{m}{k}\right)}{k^{\ell}},
\end{equation}
which we adopt as our fundamental object of study.

\begin{theorem}\label{L-thm} The differences of Hurwitz zeta values, $A_n(m,k)$
defined by~\eqref{defA}, satisfy the estimate
\begin{equation}\label{A1}
A_{n}(m,k)=\left(\frac{m}{k}-\frac{1}{2}\right)
-\frac{n}{k}\left[\psi\left(\frac{m}{k}\right)+\ln k+1-H_{n-1}\right]
+a_{n}(m,k)
\end{equation}
where the $a_n(m,k)$ are exponentially small:
\begin{equation}\label{A2}
\begin{array}{lll}
a_{n}(m,k)&=&\ds \frac{1}{k}\left(\frac{2n}{\pi k}\right)^{1/4}
\exp\left(-\sqrt{\frac{4\pi n}{k}}\right)
\cos\left(
\sqrt{\frac{4\pi n}{k}} -\frac{5\pi}{8} -\frac{2\pi m}{k} \right)
\\
&& {} \ds +\mathcal{O}\left(n^{-1/4}e^{-2\sqrt{\pi n/k}}\right).
\end{array}
\end{equation}
Here $\psi(x)=\Gamma'(x)/\Gamma(x)$ is the logarithmic derivative of the Gamma function.

\end{theorem}

The previous results for Riemann zeta may be regained by setting 
$m=k=1$, so that $\delta_n=A_n(1,1)$ and $b_n=a_n(1,1)$.

\begin{proof}
 Converting the sum to the N\"orlund-Rice integral, and extending the
contour to infinity, like before, one obtains 
\begin{equation}
A_{n}(m,k)=\frac{(-1)^{n}}{2\pi i}\, n!\,
\int_{\frac{3}{2}-i\infty}^{\frac{3}{2}+i\infty}
\frac{\zeta\left(s,\frac{m}{k}\right)}{k^{s}s(s-1)\cdots(s-n)}\, ds
\end{equation}
 Moving the contour to the left, one encounters a single pole at $s=0$
and a double pole at $s=1$. The residue of the pole at $s=0$ is
\begin{equation*}
\mbox{Res}(s=0)=\zeta\left(0,\frac{m}{k}\right)=\frac{1}{2}-\frac{m}{k}.
\end{equation*}
(See~\cite[p.~271]{WhWa27} for this evaluation.)
The double pole at $s=1$ evaluates to
\begin{equation*}
\mbox{Res}(s=1)=\frac{n}{k}
\left[\psi\left(\frac{m}{k}\right)+\ln k+1-H_{n-1}\right]
\end{equation*}
Combining these, one obtains~\eqref{A1}
where the $a_n$ are given by
\begin{equation}\label{defa}
a_{n}(m,k)=\frac{(-1)^{n}}{2\pi i}\, n!\,
\int_{-\frac{1}{2}-i\infty}^{-\frac{1}{2}+i\infty}
\frac{\zeta\left(s,\frac{m}{k}\right)}{k^{s}s(s-1)\cdots(s-n)}\, ds.
\end{equation}

As before,
the $a_n(m,k)$ have the remarkable property of being exponentially
small; that is, 
$ a_{n}(m,k)=\mathcal{O}\left(e^{-K\sqrt{n}}\right)$,
 for a constant $K$ that only depends on~$k$. 
The precise behavior of the exponentially small term may be obtained 
by using a saddle-point analysis parallel to the one given in the previous
sections. Its application here is abbreviated, as there
are no major differences in the course of the derivations.

The term $a_{n}(m,k)$ is represented by the integral of~\eqref{defa}.
At this point, the functional equation for the Hurwitz zeta may be
applied. This equation is
\begin{equation}
\zeta\left(1-s,\frac{m}{k}\right)
=\frac{2\Gamma(s)}{(2\pi k)^{s}}
\sum_{p=1}^{k}\cos\left(\frac{\pi s}{2}-\frac{2\pi pm}{k}\right)
\zeta\left(s,\frac{p}{k}\right),
\end{equation}
and it either follows from adapting the ``second proof'' of Riemann
for the common zeta function~\cite[\S2.4]{Titchmarsh86}
or from the transformation formula of Lerch's transcendent $\phi$ found in~\cite[p.~280]{WhWa27}.
This allows the integral in~\eqref{defa} to be expressed as a sum:
\begin{equation*}
a_{n}(m,k)=-\frac{2n!}{k\pi i}\sum_{p=1}^{k}
\int_{\frac{3}{2}-i\infty}^{\frac{3}{2}+i\infty}
\frac{1}{(2\pi)^{s}}\frac{\Gamma(s)\Gamma(s-1)}{\Gamma(s+n)}
\cos\left(\frac{\pi s}{2}-\frac{2\pi pm}{k}\right)
\zeta\left(s,\frac{p}{k}\right)\, ds
\end{equation*}
 It proves  convenient to pull the phase factor out of the
cosine part 
and write the integral as 
\begin{equation*}
\begin{array}{lll}
a_{n}(m,k) & = & \ds -\frac{n!}{k\pi i}
\sum_{p=1}^{k}\exp\left(i\frac{2\pi pm}{k}\right) \\
 & & \quad \ds
\int_{\frac{3}{2}-i\infty}^{\frac{3}{2}+i\infty}
\frac{1}{(2\pi)^{s}}\frac{\Gamma(s)\Gamma(s-1)}{\Gamma(s+n)}\,
\exp\left(-i\frac{\pi s}{2}\right)
\zeta\left(s,\frac{p}{k}\right)\, ds   +{\bf c.c.},
\end{array}
\label{eq:two-integrals}
\end{equation*}
 where ${\bf c.c.}$ (``\emph{complex conjugate}'') means that $i$ should be replaced by $-i$ in the
two exp parts.

%

To recast the equation above into the form needed
for the saddle-point method, an asymptotic expansion of the
integrands needs to be made for large $n$. 
As before, the appropriate scaling parameter is $x=s/\sqrt{n}$.
The asymptotic
expansion is then performed by holding $x$ constant, and taking $n$
large. Thus, one writes 
\begin{equation}\label{eq:saddle}
a_{n}(m,k) = -\frac{1}{k\pi i}\sum_{p=1}^{k}
\left[e^{i {2\pi pm}/{k}}
\int_{\sigma_{0}-i\infty}^{\sigma_{0}+i\infty}e^{\omega(x\sqrt{n})}dx\right.
\left.+e^{-i{2\pi pm}/{k}}
\int_{\sigma_{0}-i\infty}^{\sigma_{0}+i\infty}
e^{\overline{\omega}(x\sqrt{n})}dx\right].
\end{equation}

Proceeding, one finds 
\begin{equation*}
\omega(s)=\log n!+\frac{1}{2}\log n
-s\log\left(\frac{2\pi p}{k}\right)-i\frac{\pi s}{2}
+\log\frac{\Gamma(s)\Gamma(s-1)}{\Gamma(s+n)}
+\mathcal{O}\left( \left(\frac{p}{k+p}\right)^s \right)
\end{equation*}
 where the approximation
$\log \zeta\left(s,p/k\right)=(k/p)^{s}+\mathcal{O}\left( (p/(k+p))^s \right)$,
for large $\Re(s)$, has been made. 
Expanding to $\mathcal{O}(1/\sqrt{n})$
and collecting terms, one obtains 
\begin{equation}\begin{array}{lll}
\omega(x\sqrt{n}) & = & \ds \frac{1}{2}\log n
-x\sqrt{n}\left[\log\frac{2\pi p}{k}+i\frac{\pi}{2}+2-2\log x\right] \\
 &  & \ds  {} +\log2\pi-2\log x-\frac{x^{2}}{2}+\mathcal{O}\left(n^{-1/2}\right) 
\end{array}
\end{equation}

 The saddle-point is obtained  by solving $\omega'(x\sqrt{n})=0$.  To
 lowest order, one  has $x_{0}=(1+i)\sqrt{\pi  p/k}$. 
Also,  $\omega''(x\sqrt{n})=2/x\sqrt{n}+\mathcal{O}(n^{-1})$.
 Substituting in the saddle-point formula~\eqref{sadfor}, one directly finds
\begin{equation}
\int_{\sigma_{0}-i\infty}^{\sigma_{0}+i\infty}
\!\!\! e^{\omega(x\sqrt{n})}dx 
=\left(\frac{2\pi^{3}pn}{k}\right)^{1/4}\!\!\!e^{i\pi/8}
\exp\!\left(-(1+i)\sqrt{\frac{4\pi pn}{k}}\right)\! {}
+\mathcal{O}\left(n^{-1/4}e^{-2\sqrt{\pi pn/k}}\right)
\end{equation}
 while the integral for $\overline{\omega}$ is the complex conjugate 
quantity (having a saddle-point at the complex conjugate location). Inserting
this into equation \eqref{eq:saddle} gives a sum of contributions for $p=1,\ldots,k$,
of which, 
for large $n$, only the $p=1$ term is seen to contribute significantly.
So, one has the estimation~\eqref{A2} of the statement.
\end{proof}

\section{Perspective}\label{concl-sec}
\smallskip
The previous methods serve to
unify and make precise estimates carried out in the literature by a diversity
of approaches. For instance, the study of quantities arising in connection with Li's criterion
calls for estimating,  in the notations of~\eqref{defA},
\begin{equation}\label{coffeyS1}
A_n(1,2)=\sum_{\ell=2}^n \binom{n}{\ell}(-1)^\ell(1-2^{-\ell})\zeta(\ell).
\end{equation}
Bombieri and Lagarias encountered this quantity in~\cite[Th.~2]{BoLa99}
and Coffey (see his~$S_1(n)$ in~\cite{Coffey05}) proved, by means of series rearrangements akin to~\eqref{simpdel}
used in conjunction with Euler-Maclaurin summation the inequality
\begin{equation}\label{coffeyineq}
A_n(1,2) \ge \frac{n}{2}\log n+(\gamma-1)\frac{n}{2}+\frac12.
\end{equation}
Our analysis quantifies  
$A_n(1,2)$ to be
\[
A_n(1,2)=\frac{n}{2}\psi(n)+n(\gamma-\frac12+\frac12\log2)+o(1),
\]
where the $o(1)$ error term 
above is $a_n(1,2)$, which is
exponentially small and oscillating:
\begin{equation}
 a_{n}(1,2)=\frac{1}{2}\left(\frac{n}{\pi}\right)^{1/4}
\exp\left(-\sqrt{2\pi n}\right) \cos\left(\sqrt{2\pi n}-\frac{5\pi}{8}\right)
+\mathcal{O}\left(n^{-1/4}e^{-\sqrt{2\pi n}}\right)
\end{equation}

Another observation is that the combination of N\"orlund-Rice integrals
and saddle-point estimates applies to many ``desingularized'' versions of 
the Riemann zeta function, like
\[
(1-2^{1-s})\zeta(s), \quad
(s-1)\zeta(s), \quad
\zeta(2s)-\frac{1}{2s-1},\quad
(2s-1)\zeta(2s).
\]
The first one is directly amenable to Theorem~\ref{L-thm}. The
Newton series involving $\zeta(2s)$ include Ma{\'s}lanka's expansion~\cite{Maslanka01}
(relative to $(2s-1)\zeta(2s)$)
and have a striking feature---their Newton coefficients are 
polynomials in~$\pi$ with rational coefficients. 
In addition, the exponential smallness of error terms in 
asymptotic expansions of finite differences of this sort
has the peculiar feature of inducing near-identities that relate 
rational combinations of zeta values and Euler's constant.
For instance, defining the following elementary variant of~$b_{n}$,
\[
c_n:=-\sum_{k=1}^n \binom{n}{k} (-1)^k \frac{\zeta(k+1)}{k+1},
\]
we find to more than 35 digits of accuracy,
\begin{equation*}
\left\{\begin{array}{rll}
c_{499}-H_{499}+1&=&0.57\underline{\bf 8}21\,56649\,01532\,86060\,65120\,90082\,40243\ldots \\
\gamma&=&0.57\underline{\bf7}21\,56649\,01532\,86060\,65120\,90082\,40243\ldots ,
\end{array}
\right.
\end{equation*}
where the   sole discrepancy observed  is  in  the third decimal  digit.


\smallskip

The N\"orlund integrals are also of interest in the context of 
differences of inverse zeta values, for which curious relations
with the Riemann hypothesis have been noticed by Flajolet and Vall\'ee~\cite{FlVa00}, and
independently by B\'aez-Duarte  in a scholarly note~\cite{Baez05}. Consider the
typical quantity
\begin{equation}\label{invz}
d_n=\sum_{k=2}^n \binom{n}{k}(-1)^k \frac{1}{\zeta(k)},
\end{equation}
which arises as coefficient in the Newton series representation of $1/\zeta(s)$. 
Its asymptotic analysis can be approached by means of a N\"orlund-Rice representation
as noted 
in related contexts by the authors of~\cite{FlVa00} and more recently by Ma\'slanka in~\cite{Maslanka06}.
The following developments
provide a rigorous basis for some of the observations made in~\cite{Maslanka06},
simplifies the criterion for RH that is implicit in~\cite{FlVa00}, and
offers an alternative (based on real extrapolation and the Mellin transform) to B\'aez-Duarte's treatment,
while pointing in the direction of easy generalizations relative  to
$1/\zeta(2s)$, $\zeta'(s)/\zeta(s)$, 
	$\zeta(s-1)/\zeta(s)$, or other similar functions.
\begin{theorem}\label{invz-thm} 
The differences of inverse zeta values $d_n$ defined by~\eqref{invz} are such that
the following two assertions are equivalent: 
\begin{itemize}
\item[]{\bf FVBD  Hypothesis (akin to~\cite{Baez05,FlVa00})}. For any~$\epsilon>0$, there exists a constant $C_\epsilon>0$ such that
\[
|d_n|<C_\epsilon n^{1/2+\epsilon}.
\]
\item[]{\bf RH (Riemann hypothesis)}. The Riemann 
zeta function~$\zeta(s)$ is free of zeros in the half-plane~$\Re(s)>\frac12$.
\end{itemize}
\end{theorem}
\begin{proof}
$(i)$~Assume {\bf RH}. Under RH, it is known that, given any~$\sigma_0>\frac12$ and any~$\epsilon>0$, 
one has 
\begin{equation}\label{RHeps}
\frac{1}{\zeta(s)}=\O(|t|^\epsilon),
\qquad\hbox{for $\Re(s)=\sigma_0$,~~where~~$t=\Im(s)$}
\end{equation}
(see Equation~(14.2.6) of~\cite[p.~337]{Titchmarsh86}). 
Then, start from the N\"orlund integral representation
(cf Lemma~\ref{norlem} and Equation~\eqref{norbn}), 
\begin{equation}\label{norinv}
d_n=J_n(c),\quad\hbox{where}~~ J_n(c):=\frac{(-1)^{n-1}}{2\pi i}\int_{c-i\infty}^{c+i\infty}
\frac{1}{\zeta(s)} \frac{n!}{s(s-1)\cdots(s-n)}\, ds,
\end{equation}
which is valid unconditionally for $c\in(1,2)$.
Next, we propose to move the line of integration to $c=\sigma_0$.
To this effect, observe that 
the integral  $J_n(\sigma_0)$ defined in~\eqref{norinv}
converges and is $\mathcal{O}(n^{\sigma_0})$, since
\[\renewcommand{\arraystretch}{1.5}
\begin{array}{lll}
\ds\left| J_n(\sigma_0)\right| & \le & \ds
\frac{1}{2\pi\sigma_0}\binom{n-\sigma_0}{-\sigma_0}^{-1} \int_{-\infty}^\infty |\zeta(\sigma_0+it) |^{-1}\,
\left| \frac{-\sigma_0\cdots(-\sigma_0+n)}{(-\sigma_0-it)\cdots (-\sigma_0-it+n)}\right|\, dt\\
& \le & \ds
\frac{1}{2\pi\sigma_0}\binom{n-\sigma_0}{-\sigma_0}^{-1} \int_{-\infty}^\infty |\zeta(\sigma_0+it) |^{-1}\,
\left| \frac{-\sigma_0(-\sigma_0+1)}{(-\sigma_0-it)(-\sigma_0-it+1)}\right|\, dt\\
&=&\ds \O\left(n^{\sigma_0}\right).
\end{array}
\]
There, the second line results from the fact that, for~$x,t$ real, one has $|(x/(x-it)|\le1$;
the third line summarizes the asymptotic estimate  $\binom{n-\sigma_0}{-\sigma_0}^{-1}=\O(n^{\sigma_0})$ (by
Stirling's formula) as well as the fact that the integral factor is convergent (since the integrand decays at least
as fast as $\mathcal{O}(|t|^{-2+\epsilon})$ as $|t|\to+\infty$).

\smallskip

$(ii)$~Assume {\bf FVBD}. First, a reorganization similar to the one leading to~\eqref{simpdel}
but based on the expansion of $1/\zeta(s)$ shows that
\[
d_n=\sum_{\ell=1}^\infty \mu(\ell)\left[\left(1-\frac{1}{\ell}\right)^n-1+\frac{n}{\ell}\right],
\]
with $\mu(\ell)$ the M\"obius function. The general term of the sum decreases like $n/\ell^2$,
which ensures absolute convergence. Next, introduce the function
\[
D(x)=\sum_{\ell\ge1} \mu(\ell)\left[e^{-x/\ell}-1+\frac{x}{\ell}\right],
\]
whose general term decreases like $x^2/\ell^2$. 

Fix any small~$\delta>0$ ($\delta=\frac1{10}$ is suitable) and define $\ell_0=\lfloor x^{1-\delta}\rfloor$.
The difference $d_n-D(n)$ satisfies
\begin{equation}\label{dede}
\begin{array}{lll}
 d_n-D(n) &=& \ds \sum_{\ell=1}^\infty \mu(\ell)\left[\left(1-\frac{1}{\ell}\right)^n
-e^{-n/\ell}\right]
\\
&=& \ds \left(\sum_{\ell<\ell_0}+\sum_{\ell\ge \ell_0}\right)
\mu(\ell)e^{-n/\ell}\left[e^{n/\ell+n\log(1-1/\ell)}
-1\right]
\\
&=&\ds \O\left(\ell_0 e^{-n/\ell_0}\right) +\sum_{\ell\ge\ell_0}\O\left(\frac{n}{\ell^2}\right)
 \quad = \quad\ds \O\left(n^\delta\right),
\end{array}
\end{equation}
by series reorganization, a split of the sum according to $\ell \gtreqless\ell_0$,
and trivial majorizations. 

Given~\eqref{dede}, the FVBD Hypothesis implies that $D(x)=\O(x^{1/2+\epsilon})$,
at least when $x$ is a positive \emph{integer}. To extend this estimate to real values of~$x$, 
it suffices to note that $D(x)$ is differentiable on~$\R_{>0}$, and
\[
D'(x)=-\sum_{\ell=1}^\infty \frac{\mu(\ell)}{\ell}\left[e^{-x/\ell}-1\right],
\]
is proved to be $\O(1)$ by bounding techniques similar to~\eqref{dede}.
%
Thus, assuming the FVBD Hypothesis, the estimate
\begin{equation}\label{estimd}
D(x)=\O\left(x^{1/2+\epsilon}\right), \qquad x\to+\infty
\end{equation}
holds for \emph{real} values of~$x$.

Regarding the behaviour of~$D(x)$ at~$0$,
the general term of $D(x)$ is asymptotic to $x^2/(2\ell^2)$, so that $D(x)=\O(x^2)$,
as $x\to0^{+}$. This, combined with the estimate of $D(x)$ at infinity
expressed by~\eqref{estimd}, implies 
(under the FVBD Hypothesis, still) that
the Mellin transform
\begin{equation}\label{mel1}
D^\star(s):=\int_0^\infty D(x) x^{s-1}\, dx,
\end{equation}
exists and is an analytic function of $s$ for all $s$ in the strip $-2<\Re(s)<-\frac12-\epsilon$.
On the other hand, the usual properties of Mellin transforms (see, e.g., the survey~\cite{FlGoDu95})
imply that 
\begin{equation}\label{mel2}
D^\star(s)=\left(\sum_{\ell=1}^\infty \mu(\ell)\ell^{s}\right)
\cdot \int_0^\infty \left[e^{-x}-1+x\right] x^{s-1}\, dx =\frac{\Gamma(s)}{\zeta(-s)},
\end{equation}
at least for $s$ such that $-2<\Re(s)<-1$, which ensures that the 
usual series expansion of $1/\zeta(-s)$
is absolutely convergent. The comparison of the analytic character of~\eqref{mel1} in
$-2<\Re(s)<-\frac12-\epsilon$ (implied by the FVBD Hypothesis) and of the explicit form of~\eqref{mel2}
shows that the Riemann Hypothesis is a consequence of the FVBD Hypothesis.
\end{proof}

Numerically, for comparatively low values of~$n$, it would seem that 
$d_n$ tends slowly but steadily to $2=-1/\zeta(0)$. 
For instance, we have $d_{20}\doteq1.93 $, $d_{50}\doteq1.987$, $d_{100}=1.996$,
$d_{200}\doteq 1.9991$.
However, it appears from our previous analysis 
and a residue calculation applied to~\eqref{norinv} 
that there must be complicated oscillations
due to the nontrivial zeta zeros---these oscillations in fact
 eventually \emph{dominate}, though at a rather late stage,
as we now explain following~\cite{FlVa00}. 
Indeed, assuming for notational convenience the simplicity 
of the nontrivial zeros of $\zeta(s)$, one has (unconditionally)
\begin{equation}\label{inv2}
d_n= \sum_{\rho}\!{}^{\hbox{$\star$}} \frac{1}{\zeta'(\rho)} 
\frac{\Gamma(n+1)\Gamma(-\rho)}{\Gamma(n+1-\rho)}+2+o(1),
\end{equation}
where   the   summation   extends  to   all   nontrivial  zeros~$\rho$
of~$\zeta(s)$   with   $0<\Re(\rho)<1$,   while   the    starred   sum
($\sum{}\!{}^{\star}$) means  that  zeros should  be  suitably grouped,
following       the  careful  discussion     in~\S9.8     of        Titchmarsh's
treatise~\cite[p.~219]{Titchmarsh86} in  relation to  a  formula of
Ramanujan. A simplified
 model of the sequence~$d_n$ then follows from the fact that,
for large~$n$, any individual term of the sum in~\eqref{inv2} 
corresponding to a zeta zero $\rho=\sigma+i\tau$
is  asymptotically
\begin{equation}\label{zeterm}
\frac{\Gamma(-\rho)}{\zeta'(\rho)} n^{\sigma} e^{i\tau\log n}.
\end{equation}
Such a term involves  a logarithmically oscillating  component,  a
slowly growing  component  $n^\sigma$ ($\sqrt{n}$ under  RH), as well as a
multiplier that  is likely   to  be extremely small numerically, since    it involves
the quantity $\Gamma(-\rho)\asymp e^{-\pi|\tau|/2}$.  
For the first nontrivial zeta zero
at $\rho\doteq\frac12+i\,14.13$, the term~\eqref{zeterm} is very roughly
\begin{equation}\label{approxi}
10^{-9} \sqrt{n}\cos\left(14.13\, \log n\right),
\end{equation}
and for the next zero, at $\rho\doteq\frac12+i\,21.022$,
the numerical  coefficient drops to about~$10^{-14}$.
The  corresponding  oscillations
then have the curious feature of being numerically detectable only for
very large  values  of~$n$:
for instance, in order  for the first term given by~\eqref{approxi} to attain
the value~1, one needs $n\approx10^{18}$, while for the contribution
of the second zero, one would need $n\approx 10^{28}$. Since it is known that all zeros 
of~$\zeta(s)$ lie on the critical line till height $T_0\approx 5\cdot 10^8$,
we can estimate in this way that the presence of a nontrivial zero (if any) 
off the critical line could only be detected in the asymptotic behaviour 
of~$d_n$ for values of~$n$ larger than $N_0\approx e^{\pi T_0}\approx
10^{600,000,000}$.
In summary, from a numerical point of view, a possible failure  of RH,
though in theory traceable through the asymptotic behavior of $d_n$, 
is in reality violently counterbalanced by the exponential
decay of the $\Gamma$ factor, hence it must remain totally undetectable in practice---analogous    facts
  were     observed
in~\cite{FlVa00,Maslanka06}.  It is finally of  interest to note  that
such   phenomena   do  occur   in   nature,   specifically,   in   the
determination by Flajolet and Vall\'ee in~\cite{FlVa00} of  \emph{the    expected number  of  continued
fraction digits  that are  necessary to  sort $n$  real numbers
drawn uniformly at random from the unit interval}.

\medskip
\noindent
{\bf Acknowledgements.} We are thankful to Mark Coffey for several
useful comments on an original draft and to Luis B\'aez-Duarte for friendly and
informative feedback.


\def\cprime{$'$}

\end{document}